\documentclass[a4paper]{amsart}  

\usepackage{amssymb} \usepackage{amscd} \usepackage{times}

\def\zbb{\mathbb{Z}}  
  
  \def\phi{\varphi}
 \def\p1{{\mathbb{P}^1_\zbb}}

\newtheorem{Theorem}{\quad Theorem}[section]

\newtheorem{Corollary}[Theorem]{\quad Corollary} 

\newtheorem{Lemma}[Theorem]{\quad Lemma}

\begin{document}

\title{ Uniqueness type result in dimension 3.}

\author{Samy Skander Bahoura}

\address{Departement de Mathematiques, Universite Pierre et Marie Curie, 2 place Jussieu, 75005, Paris, France.}
              
\email{samybahoura@yahoo.fr} 

\date{}

\maketitle

\begin{abstract}

We give some estimates of type $ \sup \times \inf $ on Riemannian manifold of dimension 3 for the prescribed curvature type equation. As a consequence, we derive an uniqueness type result.

\end{abstract}

\section{Introduction and Main Results} 

\bigskip

In this paper, we deal with the following prescribed scalar curvature type equation in dimension 3:

$$ \Delta u+h(x)u=V(x) u^5,\,\, u >0 . \qquad (E) $$

Where $ h, V $ are two continuous functions. In the case  $ 8 h= R_g $ the scalar curvature, we call $ V $ the prescribed scalar curvature. Here, we assume $ h $ a bounded function and $ h_0=||h||_{L^{\infty}(M)} $.

\bigskip

We consider three positive real number $ a,b, A $ and we suppose $ V $ lipschitzian:

 $$ 0 < a \leq V(x) \leq b < + \infty \,\, {\rm and} \,\, ||\nabla V||_{L^{\infty}(M)} \leq A. \qquad (C) $$

\bigskip

The equation $ (E) $ was studied a lot, when $ M =\Omega \subset {\mathbb R}^n $ or $ M={\mathbb S}_n $ see for example, [2-4], [11], [15]. In this case we have a $ \sup \times \inf $ inequality.

\smallskip

The corresponding equation in two dimensions on open set $ \Omega $ of $ {\mathbb R}^2 $, is:

$$ \Delta u=V(x)e^u, \qquad (E') $$

The equation $ (E') $ was studed by many authors and we can find very important result about a priori estimates in [8], [9], [12], [16], and [19]. In particular in [9] we have the following interior estimate:

$$  \sup_K u  \leq c=c(\inf_{\Omega} V, ||V||_{L^{\infty}(\Omega)}, \inf_{\Omega} u, K, \Omega). $$

And, precisely, in [8], [12], [16], and [19], we have:

$$ C \sup_K u + \inf_{\Omega} u \leq c=c(\inf_{\Omega} V, ||V||_{L^{\infty}(\Omega)}, K, \Omega), $$

and,

$$ \sup_K u + \inf_{\Omega} u \leq c=c(\inf_{\Omega} V, ||V||_{C^{\alpha}(\Omega)}, K, \Omega). $$

where $ K $ is a compact subset of $ \Omega $, $ C  $ is a positive constant which depends on $\dfrac{\inf_{\Omega} V}{\sup_{\Omega} V} $, and,  $ \alpha \in (0, 1] $.

\smallskip

In the case $ V\equiv 1 $ and $ M $ compact, the equation $ (E) $ is Yamabe equation. Yamabe has tried to solve problem but he could not, see [22]. T.Aubin and R.Schoen have proved the existence of solution in this case, see for example [1] and [14] for a complete and detailed summary.

\bigskip

When $ M $ is a compact Riemannian manifold, there exist some compactness result for equation  $ (E) $ see [18]. Li and Zhu see [18], proved that the energy is bounded and if we suppose $ M $ not diffeormorfic to the three sphere, the solutions are uniformly bounded. To have this result they use the positive mass theorem.

\bigskip

Now, if we suppose $ M $ Riemannian manifold (not necessarily compact) and $ V\equiv 1 $, Li and Zhang [17] proved that the product $ \sup \times \inf $ is bounded. Here we extend the result of [6].

 Our proof is an extension of  Brezis-Li and Li-Zhang result in dimension 3, see [7] and [17], and,  the moving-plane method is used to have this estimate. We refer to Gidas-Ni-Nirenberg for the  moving-plane method, see  [13]. Also, we can see in [10], one of the application of this method.
 
\bigskip

Here, we give an equality of type $ \sup \times \inf $ for the equation $ (E) $ with general conditions $ (C) $.  Note that, in our proof, we do not need a classification result for some particular elliptic PDEs on $ {\mathbb R }^3 $.

In dimension greater than 3 we have other type of estimates by using moving-plane method, see for example [3, 5].

There are other estimates of type $ \sup + \inf $ on complex Monge-Ampere equation on compact  manifolds, see [20-21] . They consider, on compact Kahler manifold $ (M, g) $, the following equation:

\begin{center}

$ \begin{cases} 

(\omega_g+\partial \bar \partial \phi)^n=e^{f-t\phi} \omega_g^n, \\

\omega_g+\partial \bar \partial \phi >0 \,\, {\rm on } \,\,  M  \\
 
\end{cases}  $ 
    
\end{center} 

And, they prove some estimates of type $ \sup_M+m \inf_M \leq C $ or $  \sup_M + m \inf_M \geq C $ under the positivity of the first Chern class of M.

\bigskip

Here, we have,

\bigskip

\begin{Theorem}  For all compact set $ K $ of $ M $ and all positive numbers $ a,b,A, h_0 $ there is a positive constant c, which depends only on, $ a, b, A, h_0, K, M, g $ such that:

$$ \sup_K  u  \times \inf_M u \leq c, $$

for all $ u $ solution of $ (E) $ with conditions $ (C) $.

\end{Theorem} 

This theorem generalise Li-Zhang result, see [17] in the case $ V\equiv 1 $. Here, we use Li and Zhang method in [17].

\bigskip

In the case $ h\equiv \epsilon \in (0,1) $ and $ u_{\epsilon} $ solution of :

$$  \Delta u_{\epsilon} +\epsilon u_{\epsilon} =V_{\epsilon} u_{\epsilon} ^5,\,\, u_{\epsilon} >0 . \qquad (E_{\epsilon} ) $$

We have:

\begin{Corollary}  For all compact set $ K $ of $ M $ and all positive numbers $ a,b,A $ there is a positive constant c, which depends only on, $ a, b, A, K, M, g $ such that:

$$ \sup_K  u_{\epsilon}  \times \inf_M u_{\epsilon} \leq c, $$

for all $ u $ solution of $ (E_{\epsilon}) $ with conditions $ (C) $.

\end{Corollary} 

Now, if we assume $ M $ a compact riemannian manifold and $ 0 < a \leq V_{\epsilon} \leq b < +  \infty $

we have:

\begin{Theorem}({\rm see [3]}). For all positive numbers $ a, b, m  $ there is a positive constant c, which depends only on, $ a, b, m, M, g $ such that:

$$ \epsilon  \sup_M  u_{\epsilon}  \times \inf_M u_{\epsilon} \geq c, $$

for all $ u_{\epsilon} $ solution of $ (E_{\epsilon}) $ with

$$  \max_M u_{\epsilon} \geq m >0 $$.

\end{Theorem} 

As a consequence of the two previous theorems, we have:

\begin{Theorem}  For all positive numbers $ a, b, A $ we have:

$$  \max_M u_{\epsilon} \to 0, $$ and (up to a subsequence),

$$ \dfrac {\max_M u_{\epsilon}}{{\epsilon}^{1/4}} \to w_0 >0, \,\, \,\, {\rm and,} \,\, \,\, \dfrac {\min_M u_{\epsilon}}{{\epsilon}^{1/4}} \to w_0 >0. $$

\end{Theorem} 
 
 \bigskip
 
 \underbar {\bf  Remarks:} 

\bigskip

\begin{itemize}

\item  It is not necessary to have $ u_{\epsilon} \equiv w_0{\epsilon}^{1/4} $, because if we take a nonconsant function $ V $, we can find by the variational approach a non constant positive solution of the subcritical equation:
$$  \Delta u_{\epsilon} + \epsilon u_{\epsilon} = {\mu}_{\epsilon}V u_{\epsilon}^{5- \epsilon} , \,\, {\rm with} \,\, {\mu}_{\epsilon}, u_{\epsilon} >0. $$
In this case (subcritical which tends to the critical) we also have the $ \sup \times \inf  $ inequalities and the uniqueness type theorem.

\item In fact, we prove, up to a subsequence that $ \dfrac {u_{\epsilon}}{{\epsilon}^{1/4}} $ converge to a constant which depends on $ a, b $ and $ A $.

\end{itemize}

\section{Proof of the theorems} 

\underbar {\it Proof of theorem 1.1:} 

\bigskip

We want to prove that:

$$ \epsilon \max_{B(0, \epsilon)} u \times \min_{B(0, 4 \epsilon) } u \leq c=c(a, b, A, M, g) \qquad (1) $$

We argue by contradiction and we assume that:

$$  \max_{B(0, \epsilon_k)} u_k \times \min_{B(0, 4 \epsilon_k) } u_k \geq k {\epsilon_k }^{-1} \qquad (2)$$

\underbar {Step 1: The blow-up analysis }

\bigskip

The blow-up analysis gives us :

For some $  \bar x_k \in B(0, \epsilon_k) $, $ u_k( \bar x_k)=\max_{B(0, \epsilon_k)} u_k $, and, from the hypothesis,

$$ u_k( \bar x_k)^2 \epsilon_k \to + \infty. $$

By a standard selection process, we can find $ x_k  \in B(\bar x_k , \epsilon_k/2) $ and $ \sigma_k \in (0, \epsilon_k/4) $ satisfying,

$$ u_k(x_k)^2\sigma_k \to + \infty, \qquad (3)$$

$$ u_k(x_k) \geq  u_k(\bar x_k), \qquad (4)$$

and,

$$ u_k(x) \leq  C_1 u_k (x_k),  \,\, {\rm in} \,\, B(x_k, \sigma_k), \qquad (5)$$

where $ C_1 $ is some universal constant.

It follows from above ( $ (2), (4) $) that:

$$ u_k(x_k) \times \min_{\partial B(x_k, 2 \epsilon_k) } u_k \epsilon_k  \geq  u_k(\bar x_k) \times \min_{B(0, 4 \epsilon_k) } u_k \epsilon_k \geq k   \to + \infty .\qquad (6) $$

We use $ \{ z^1, \ldots, z^n \}  $ to denote some geodesic normal coordinates centered at $ x_k $ (we use the exponential map). In the geodesic normal coordinates, $ g=g_{ij}(z)dz"dz^j $,

$$ g_{ij}(z)={\delta}_{ij} = O(r^2), \,\, g:=det(g_{ij}(z))=1+O(r^2), \,\, h(z)=O(1), \qquad (7)$$

where $ r= |z| $. Thus,

$$ \Delta_g u=\dfrac{1}{\sqrt g} \partial_i ( \sqrt g g^{ij} \partial_j u)=\Delta u+ b_i\partial_i u +  d_{ij} \partial_{ij} u, $$

where

$$ b_j=O(r),  \,\, \partial_{ij}= O(r^2) \qquad (8)$$

We have a new function:

$$ v_k(y)=M_k^{-1} u_k(M_k^{-2}y) \,\, {\rm for} \,\, |y| \leq 3 \epsilon_k  M_k^{2} $$

where $ M_k=u_k(0) $.

From $ (5), (6) $, we have:

\bigskip

 $ \begin{cases} 

 \Delta v_k + \bar b_i\partial_i v_k + \bar d_{ij} \partial_{ij} v_k - \bar c v_k +{v_k}^{5} =0  \,\, {\rm for} \,\, |y| \leq 3 \epsilon_k M_k^{2}  \\
 
 v_k(0)=1  \\
  
  v_k(y) \leq C_1 \, {\rm for} \,\, |y| \leq  \sigma_k  M_k^{2} \\
  
 \lim_{k \to + \infty } \min_{|y|=2\epsilon_k M_k^{2}} (v_k(y)|y| )= + \infty  \qquad (9) \\

 \end{cases}  $

where $ C_1 $ is a universal constant and

$$ \bar b_i(y)= M_k^{-2} b_i(M_k^{-2}y), \,\, \bar d_{ij}(y)=d_{ij}(M_k^{-2}y) \qquad (10)$$ 

and,

$$ \bar c(y)=M_k^{-4} h(M_k^{-2}y) \qquad (11)$$

We can see that for  $ |y| \leq 3 \epsilon_k M_k^{2}  $, we have:

$$ |\bar b_i(y)| \leq CM_k^{-4}|y|, \,\, |\bar  d_{ij}(y)| \leq CM_k^{-4} |y|^2, \,\, |\bar c(y)| \leq CM_k^{-4} \qquad (12) $$

where $ C $ depends on $ n, M, g $.

\bigskip

It follows from $ (9), (10), (11), (12) $ and the elliptic estimates, that, along a subsequence, $ v_k $ converges in $ C^2 $ norm on any compact subset of $ {\mathbb R}^2 $ to a positive function $ U $ satisfying:

\bigskip

 \begin{center}

$ \begin{cases} 

\Delta U+ U^{5} =0  \,\, {\rm in } \,\,  {\mathbb R}^2  \\
 
U(0)=1, \,\, 0 < U \leq  C_1    \qquad (13)  
\end{cases}  $ 
    
\end{center} 

\underbar {Step 2: The Kelvin transform and moving-plane method}

For $ x\in {\mathbb R}^2 $ and $ \lambda >0 $, let,

$$ v_k^{\lambda, x}(y):= \dfrac {\lambda}{|y-x|}v_k \left (x+\dfrac {\lambda^2(y-x)}{|y-x|^2}\right  ) $$

denote the Kelvin transformation of $ v_k $ with respect to the ball centered at x and of radius $  \lambda  $.

\bigskip

We want to compare for fixed $ x $, $ v_k $ and $ v_k^{\lambda, x} $. For simplicity we assume $ x=0 $. We have:

$$ v_k^{\lambda}(y):= \dfrac {\lambda}{|y|}v_k(y^{\lambda}), \,\, {\rm with } \,\,  y^{\lambda}= \dfrac {\lambda^2y}{|y|^2} $$

For $ \lambda >0 $, we set,

$$ \Sigma_{\lambda}=B\left (0, \epsilon_k{M_k}^{2}\right )-{\bar B(0,\lambda)}. $$

The boundary condition, $ (9) $, become:

$$ \lim_{k \to + \infty } \min_{|y|=\epsilon_k M_k^{2}} \left( v_k(y)|y| \right )= \lim_{k \to + \infty } \min_{|y|=2\epsilon_k M_k^{2}} \left ( v_k(y)|y| \right )= + \infty \qquad (14) $$

We have:

$$ \Delta v_k^{\lambda} + V_k^{\lambda} (v_k^{\lambda})^5 = E_1(y) \,\, {\rm for  } \,\, y \in \Sigma_{\lambda } \qquad (15) $$
 
where,

$$ E_1(y)=- \left (\dfrac{\lambda}{|y|} \right )^5 \left (  \bar b_i(y^{\lambda}) \partial_i v_k(y^{\lambda}) + \bar d_{ij}(y^{\lambda}) \partial_{ij} v_k(y^{\lambda}) - \bar c(y^{\lambda}) v_k(y^{\lambda}) \right ). \qquad (16)  $$

Clearly, from $ (10), (11) $, there exists $ C_2= C_2(\lambda_1 ) $ such that,

$$ | E_1(y) |  \leq C_2 \lambda^5 M_k^{-4} |y|^{-5} \,\, {\rm for  } \,\, y \in \Sigma_{\lambda } \qquad (17) $$

Let,

$$ w_{\lambda} = v_k- v_k^{\lambda }. $$

Here, we have, for simplicity, omitted $ k $. We observe that by $ (9), (15) $:

$$ \Delta w_{\lambda} + \bar b_i\partial_i w_{\lambda} + \bar d_{ij} \partial_{ij} w_{\lambda}- \bar c w_{\lambda} + 5\xi^4 V_kw_{\lambda} = E_{\lambda} \,\, {\rm in  }  \,\, \,\,    \Sigma_{\lambda } \qquad (18) $$

where $ \xi $ stay between $ v_k $ and $ v_k^{\lambda} $, and,

$$ E_{\lambda} = -  \bar b_i \partial_i v_k^{\lambda} + \bar d_{ij} \partial_{ij} v_k^{\lambda} + \bar c v_k^{\lambda} -E_1-(V_k-V_k^{\lambda})(v_k^{\lambda})^5  . \qquad (19)$$

A computations give us the following two estimates:

$$ |\partial_i v_k^{\lambda}(y)|\leq C \lambda |y|^{-2},  \,\, {\rm and   }  \,\,\,\, |\partial_{ij} v_k^{\lambda}(y)|  \leq C \lambda |y|^{-3} \,\, {\rm in  }  \,\, \Sigma_{\lambda } \qquad (20) $$

From $ (10), (11), (20) $, we have,

\begin{Lemma}. For somme constant $ C_3 = C_3 (\lambda ) $

$$  | E_{\lambda} | \leq C_3 \lambda M_k^{-4} |y|^{-1} +C_3 \lambda^5 M_k^{-2} |y|^{-4}  \,\, {\rm in  }  \,\, \,\, \Sigma_{\lambda }  \qquad (21) $$

\end{Lemma} 

we consider the following auxiliary function:

$$ h_{\lambda} = -C_1AM_k^{-2} {\lambda}^2 \left ( 1- \dfrac {\lambda}{|y|} \right )+C_2AM_k^{-2}{\lambda}^3\left (1- \left (\dfrac {\lambda}{|y|}\right )^2\right )-C_3M_k^{-4}\lambda(|y| - \lambda), $$

where $ C_1, C_2 $ and $ C_3 $ are three positive numbers.

\begin{Lemma}. We have,

$$ w_{\lambda } + h_{\lambda } \geq 0 , \,\, {\rm in  }  \,\, \Sigma_{\lambda }  \,\, \forall 0 < \lambda \leq \lambda_1 \qquad (22) $$
 
\end{Lemma} 

Proof of Lemma 2.2. We divide the proof into two steps.

\bigskip

\underbar {Step 1}. There exists $ \lambda_{0, k } >0 $ such that $ (22) $ holds :

$$ w_{\lambda } + h_{\lambda } \geq 0 , \,\, {\rm in  }  \,\, \Sigma_{\lambda }  \,\, \forall 0 < \lambda \leq \lambda_{0, k }. $$

To see this, we write:

$$ w_{\lambda }= v_k(y)- v_k^{\lambda }(y)=\dfrac{1}{\sqrt |y|} \left (\sqrt |y| v_k(y)- \sqrt |y^{\lambda }|v_k(y^{\lambda })\right ). $$

Note that $ y $ and $ y^{\lambda } $ are on the same ray starting from the origin. Let, in polar coordinates,

$$ f(r,\theta) =\sqrt r v_k(r,\theta). $$

From the uniform convergence of $ v_k $, there exists $ r_0 >0 $ and $ C >0 $ independant of $ k $ such that,

$$ \dfrac{\partial f}{\partial r}(r,\theta) > Cr^{-1/2} \,\, {\rm for  }  \,\,  0 < r < r_0. $$

Consequently, for $ 0 < \lambda  < |y| <r_0 $, we have:

$$ w_{\lambda }(y) + h_{\lambda } (y)=v_k(y)- v_k^{\lambda }(y)+ h_{\lambda } (y), $$

$$ >\dfrac{1}{\sqrt r_0} C{\sqrt r_0}^{-1/2 }(|y| - |y^{\lambda }|)+h_{\lambda } (y) $$

$$ >(\dfrac{C}{r_0} -C_3\lambda M_k^{-2})(|y| - \lambda ) \,\, {\rm since  }  \,\,  |y| - |y^{\lambda }| > |y| -\lambda $$

$$ >0. \qquad (23) $$

Since,

$$ |h_{\lambda }(y) | +v_k^{\lambda }(y) \leq C(k,r_0) \lambda, \,\,\,   r_0\leq |y| \leq \epsilon_k {M_k }^{2 }, $$

we can pick small $ \lambda_{0, k } \in (0, r_0) $ such that for all $ 0 < \lambda \leq \lambda_{0, k } $ we have,

$$ w_{\lambda }(y) + h_{\lambda } (y) \geq \min_{B(0, \epsilon_kM_k^{2} ) } v_k- C(k, r_0) \lambda_{0,k} >0  \,\,  \forall  \,\, r_0\leq |y| \leq \epsilon_k {M_k }^{2 } $$

Step 1 follows from $ (23) $.

\bigskip

Let,

$$ \bar {\lambda }^k = \sup \{ 0  < \lambda \leq \lambda_1, w_{\mu} + h_{\mu} \geq 0 , \,\, {\rm in  }  \,\, \Sigma_{\mu }  \,\, \forall 0 < \mu \leq \lambda \} \qquad (24) $$

\underbar {Step 2}. $ \bar {\lambda }^k = \lambda_1 $, $ (22) $ holds.

For this, the main estimate needed is:

$$  (\Delta + \bar b_i\partial_i  + \bar d_{ij} \partial_{ij} - \bar c  + 5\xi^4 V_k) ( w_{\lambda} + h_{\lambda })\leq 0 \,\, {\rm in  }  \,\, \,\,    \Sigma_{\lambda } \qquad (25) $$

Thus,

$$ \Delta h_{\lambda} + \bar b_i\partial_i h_{\lambda} + \bar d_{ij} \partial_{ij} h_{\lambda}+ (- \bar c + 5\xi^4 V_k)h_{\lambda}+ E_{\lambda }\leq 0 \,\, {\rm in  }  \,\, \,\,    \Sigma_{\lambda } . \qquad (26) $$

We have $ h_{\lambda} < 0 $, and, $ (12) $  and a computation give us,

$$ | \bar c h_{\lambda}| \leq C_3 \lambda M_k^{-4} |y|^{-1} + C_3 \lambda^2 M_k^{-6} \leq C_3 \lambda M_k^{-4} |y|^{-1} , $$

and,

$$ |\bar b_i\partial_i h_{\lambda}| + |\bar d_{ij} \partial_{ij} h_{\lambda}| \leq C_3 \lambda M_k^{-8} |y| + C_3 \lambda^3 M_k^{-6} |y|^{-1} + C_3 \lambda^5 M_k^{-6} |y|^{-2}, $$

$$ \leq C_3 \lambda M_k^{-4} |y|^{-1} +C_3 \lambda^5 M_k^{-2} |y|^{-4}  $$

Thus,

$$ |\bar b_i\partial_i h_{\lambda}| + |\bar d_{ij} \partial_{ij} h_{\lambda}|+| \bar c h_{\lambda}| \leq C_3 \lambda M_k^{-4} |y|^{-1} +C_3 \lambda^5 M_k^{-2} |y|^{-4}  \,\, {\rm in  }  \,\, \,\, \Sigma_{\lambda }  $$

Thus, by $ (21) $, 

$$ \Delta h_{\lambda} + \bar b_i\partial_i h_{\lambda} + \bar d_{ij} \partial_{ij} h_{\lambda}+ (- \bar c + 5\xi^4 V_k)h_{\lambda}+ E_{\lambda } \leq $$

$$ \leq  \Delta h_{\lambda} + C_3 \lambda M_k^{-4} |y|^{-1} + C_3 \lambda^5 M_k^{-2} |y|^{-4}+ |E_{\lambda}| \leq 0, $$

because,

$$ \Delta h_{\lambda} =-2C_3 \lambda M_k^{-4} |y|^{-1} -2C_3 \lambda^5 M_k^{-2} |y|^{-4}. $$

From the boundary condition and the definition of $ v_k^{\lambda } $ and $ h_{\lambda} $, we have:

$$  | h_{\lambda } (y)| +v_k^{\lambda }(y) \leq  \dfrac{ C(\lambda_1) }{|y| } ,   \,\,\,\, \forall   \,\,  |y|= \epsilon_k {M_k }^2,    $$

and, thus,

$$  w_{\bar \lambda^k }(y) + h_{ \bar \lambda^k } (y) >0 \,\,\,\, \forall   \,\,  |y|= \epsilon_k {M_k }^2,    $$

We can use the maximum principal and the Hopf lemma to have:

$$ w_{\bar \lambda^k } + h_{ \bar \lambda^k } >0, \,\, {\rm in  }  \,\, \,\, \Sigma_{\lambda } , $$

and,

$$  \dfrac{\partial }{ \partial \nu} (w_{\bar \lambda^k } + h_{ \bar \lambda^k }) >0, \,\, {\rm in  }  \,\, \,\, \Sigma_{\lambda }. $$

From $ (25) $ and above we conclude that $ \bar \lambda^k = \lambda_1 $ and lemma 2.2 is proved.

\bigskip

Given any  $ \lambda >0 $, since  the sequence $ v_k $ converges to $ U $ and $ h_{ \bar \lambda^k } $ converges to 0 on any compact subset of $ {\mathbb R}^2 $, we have:

$$ U(y) \geq U^{\lambda}(y), ,\,\,\, \forall   \,\,  |y| \geq \lambda, \,\, \forall   \,\,  0 < \lambda  < \lambda_1. $$

Since $ \lambda_1>0 $ is arbitrary, and since we can apply the same argument to compare $ v_k $ and $ v_k^{\lambda, x} $, we have:

$$ U(y) \geq U^{\lambda, x}(y), ,\,\,\, \forall   \,\,  |y-x| \geq \lambda >0. $$

Thus implies that $ U $ is a constant which is a contradiction.

\bigskip 

\underbar {\it Proof of theorem 1.4:} 

\bigskip 

From theorem 2.1 (see [3]), we have:

$$ \max_M u_{\epsilon} \to 0.  \qquad  (27) $$

We conclude with the aid of the elliptic estimates and the classical Harnack inequality that:

$$  \max_M u_{\epsilon}  \leq C \min_M u_{\epsilon},   \qquad (28) $$

where $ C $ is a positive constant independant of $  \epsilon $.

\bigskip

 Let $ G_{\epsilon}  $ the Green function of the operator $ \Delta+\epsilon $, we have,

$$ \int_M G_{\epsilon}(x,y)dV_g(y)=\dfrac{1}{\epsilon}, \,\, \forall \,\, x\in M. \qquad (29) $$

We write:

$$ \inf_M u_{\epsilon} =u_{\epsilon}(x_{\epsilon})=\int_M G_{\epsilon}(x_{\epsilon},y)V_{\epsilon}(y)u_{\epsilon} ^5(y)dV_g(y) \geq $$

$$  \geq a (\inf_M u_{\epsilon} )^5\int_M G_{\epsilon}(x_{\epsilon},y)dV_g(y) =a \dfrac{(\inf_M u_{\epsilon})^{5}}{\epsilon},$$

thus,

$$ \inf_M u_{\epsilon}  \leq C_1 {\epsilon}^{1/4}. \qquad (30)$$

With the similar argument, we have :

$$ \sup_M u_{\epsilon} \geq  C_2 {\epsilon}^{1/4}.\qquad (31) $$
\bigskip

Finaly, we have:

$$ C_1 {\epsilon}^{1/4}\leq u_{\epsilon} (x)  \leq C_2{\epsilon}^{1/4} \,\,\, \forall \,\, x \in  M. \qquad (32) $$

Where $ C_1 $ and $ C_2 $ are two positive constant independant of $ \epsilon $. 

\bigskip 

We set $ w_{\epsilon} =\dfrac{u_{\epsilon}}{\epsilon^{1/4}} $, then,

$$ \Delta w_{\epsilon}+ \epsilon w_{\epsilon}= \epsilon V_{\epsilon} w_{\epsilon}^5. $$

The theorem follow from the standard elliptic estimate, the Green function of the lapalcian and the Green representation formula for the solutions of the previous equation.


\bigskip 

\begin{thebibliography}{99} 

\bibitem{1}{T. Aubin. Some Nonlinear Problems in Riemannian Geometry. Springer-Verlag 1998 }


\bibitem{2}{ S.S Bahoura. Majorations du type $ \sup u \times \inf u \leq c $ pour l'\'equation de la courbure scalaire sur un ouvert de $ {\mathbb R}^n, n\geq 3 $. J. Math. Pures. Appl.(9) 83 2004 no, 9, 1109-1150.}

\bibitem{3}{S.S. Bahoura. Harnack inequalities for Yamabe type equations.  Bull. Sci. Math.  133  (2009),  no. 8, 875-892}

\bibitem{4}{S.S. Bahoura. Lower bounds for sup+inf and sup $ \times $ inf and an extension of Chen-Lin result in dimension 3.  Acta Math. Sci. Ser. B Engl. Ed.  28  (2008),  no. 4, 749-758}

\bibitem{5}{S.S. Bahoura. Estimations uniformes pour l'Žquation de Yamabe en dimensions 5 et 6. J. Funct. Anal.  242  (2007),  no. 2, 550-562.}

\bibitem{6}{S.S. Bahoura. sup $ \times $ inf  inequality on manifold of dimension 3, to appear in MATHEMATICA AETERNA}

\bibitem{7}{H. Brezis, YY. Li. Some nonlinear elliptic equations have only constant solutions. J. Partial Differential Equations 19 (2006), no. 3, 208Ð217.
}


\bibitem{7}{H. Brezis, YY. Li , I. Shafrir. A sup+inf inequality for some
nonlinear elliptic equations involving exponential
nonlinearities. J.Funct.Anal.115 (1993) 344-358.
}

\bibitem{8}{H.Brezis and F.Merle, Uniform estimates and blow-up bihavior for solutions of $ -\Delta u=Ve^u $ in two dimensions, Commun Partial Differential Equations 16 (1991), 1223-1253.
}

\bibitem{9}{L. Caffarelli, B. Gidas, J. Spruck. Asymptotic symmetry and local
behavior of semilinear elliptic equations with critical Sobolev
growth. Comm. Pure Appl. Math. 37 (1984) 369-402.
}

\bibitem{10}{C-C.Chen, C-S. Lin. Estimates of the conformal scalar curvature
equation via the method of moving planes. Comm. Pure
Appl. Math. L(1997) 0971-1017.}

\bibitem{11}{C-C.Chen, C-S. Lin. A sharp sup+inf inequality for a nonlinear elliptic equation in ${\mathbb R}^2$.
Commun. Anal. Geom. 6, No.1, 1-19 (1998).}


\bibitem{12}{B. Gidas, W-Y. Ni, L. Nirenberg. Symmetry and related properties via the maximum principle.  Comm. Math. Phys.  68  (1979), no. 3, 209-243.}


\bibitem{13}{J.M. Lee, T.H. Parker. The Yamabe problem. Bull.Amer.Math.Soc (N.S) 17 (1987), no.1, 37 -91.}

\bibitem{14}{YY. Li. Prescribing scalar curvature on $ {\mathbb S}_n $ and related
Problems. C.R. Acad. Sci. Paris 317 (1993) 159-164. Part
I: J. Differ. Equations 120 (1995) 319-410. Part II: Existence and
compactness. Comm. Pure Appl.Math.49 (1996) 541-597.
}


\bibitem{15}{YY. Li. Harnack Type Inequality: the Method of Moving Planes. Commun. Math. Phys. 200,421-444 (1999).}

\bibitem{16}{YY. Li, L. Zhang. A Harnack type inequality for the Yamabe equation in low dimensions.  Calc. Var. Partial Differential Equations  20  (2004),  no. 2, 133--151.}

\bibitem{17}{YY.Li, M. Zhu. Yamabe Type Equations On Three Dimensional Riemannian Manifolds. Commun.Contem.Mathematics, vol 1. No.1 (1999) 1-50.}

\bibitem{18}{I. Shafrir. A sup+inf inequality for the equation $ -\Delta u=Ve^u $. C. R. Acad.Sci. Paris S\'er. I Math. 315 (1992), no. 2, 159-164.}

\bibitem{19}{Y-T. Siu. The existence of Kahler-Einstein metrics on manifolds with positive anticanonical line bundle and a suitable finite symmetry group.  Ann. of Math. (2)  127  (1988),  no. 3, 585Ð627}

\bibitem{20}{G. Tian. A Harnack type inequality for certain complex Monge-Ampre equations.  J. Differential Geom.  29  (1989),  no. 3, 481Ð488.}

\bibitem{20}{Yamabe, H. On a deformation of Riemannian structures on compact manifolds. Osaka Math. J. 12 1960 21Ð37.}


\end{thebibliography}
\end{document}